\documentclass[11pt,a4paper,openany,oneside]{article}
\usepackage[a4paper,left=2.5cm,right=2.5cm, top=3cm, bottom=3cm]{geometry}
\usepackage{times}
\usepackage{authblk}
\usepackage[latin1]{inputenc}
\usepackage[russian, english]{babel}
\usepackage{upgreek}
\usepackage{mathrsfs}
\usepackage{stmaryrd}
\usepackage{amssymb,amsmath,amsthm}
\usepackage{hyperref}
\usepackage{eurosym}
\usepackage{color}
\usepackage{enumerate}
\usepackage{multirow}
\usepackage{mathtools}
\usepackage{scalerel}
\usepackage{booktabs}
\usepackage{float}
\usepackage[T2A,T1]{fontenc}
\usepackage{amsfonts}
\usepackage{cases}
\usepackage{thmtools}

\newtheorem{thm}{Theorem}[section]
\newcounter{maintheorem}

\newtheorem{mainth}[maintheorem]{Theorem}

\newtheorem{lem}[thm]{Lemma}
\newtheorem{lemma}[thm]{Lemma}

\theoremstyle{definition}
\newtheorem{defn}{Definition}[section]
\theoremstyle{remark}
\newtheorem{remark}{Remark}
\definecolor{myred}{rgb}{0.84,0.07,0.14}

\newcommand{\R}{\mathbb R}
\newcommand{\N}{\mathbb N}

\newcommand{\les}{\lesssim}
\newcommand{\intN}{\int_{\R^N}}

\newcommand{\intOmega}{\int_\Omega}

\newcommand{\divv}{{\rm div}}

\def\e{{\rm e}}

\newcommand{\cM}{{\mathcal M}}

 \def\dd{\, {\rm d}}

\newcommand{\CyrB}{\mbox{\usefont{T2A}{\rmdefault}{m}{n}\CYRB}}

\newcommand{\tb}{{\widetilde b}}
\newcommand{\tp}{{\widetilde p}}

\newcommand{\tQ}{{\widetilde Q}}
\newcommand{\tw}{{\widetilde w}}

\DeclareOldFontCommand{\it}{\normalfont\itshape}{\mathit}

\newcommand{\Erad}{E_{\text{rad}}}
\newcommand{\ges}{\gtrsim}

\numberwithin{equation}{section}

\begin{document}
	
	\renewcommand{\thefootnote}{\fnsymbol{footnote}}
	\footnotetext{\emph{Keywords:} Choquard equation, zero mass, exponential growth, variational methods, limiting Sobolev embeddings.}
	\renewcommand{\thefootnote}{\fnsymbol{footnote}}
	\footnotetext{\emph{Mathematics Subject Classification 2020:} 35A15, 35J20, 35J60, 35B33.}
	\renewcommand{\thefootnote}{\arabic{footnote}}
	
	\title{Choquard equations with critical exponential nonlinearities\\ in the zero mass case}
	
	\date{\today}
	\author{Giulio Romani}

	\affil{\small Dipartimento di Scienza e Alta Tecnologia \protect\\ Universit\`{a} degli Studi dell'Insubria\protect\\ and\protect\\ RISM-Riemann International School of Mathematics\protect\\ Villa Toeplitz, Via G.B. Vico, 46 - 21100 Varese, Italy \protect\\\texttt{giulio.romani@uninsubria.it}}
	
	\maketitle

	
	\begin{abstract}{We investigate Choquard equations in $\R^N$ driven by a weighted $N$-Laplace operator and with polynomial kernel and zero mass. Since the setting is limiting for the Sobolev embedding, we work with nonlinearities which may grow up to the critical exponential. We establish existence of a positive solution by variational methods, completing the analysis in \cite{R}, where the case of a logarithmic kernel was considered.}
	\end{abstract}
	
	\section{Introduction}
	Aim of this work is to study the weighted Choquard equation with zero mass and polynomial kernel given by
	\begin{equation}\label{Choq_mu}\tag{$\text{Ch}_0$}
		-\divv\left(A(|x|)|\nabla u|^{N-2}\nabla u\right)=\left(\frac1{|\cdot|^\mu}\ast Q(|\cdot|)F(u)\right)Q(|x|)\,f(u)\quad\ \mbox{in}\ \ \R^N.
	\end{equation}
	with $N\geq2$. Here $A$ and $Q$ are positive radial weight functions, $\mu\in(0,N)$, and the nonlinearity is positive. Since the operator is built on the $N$-Laplacian, one expects that the maximal integrability for the nonlinearity $f$ is exponential. This is indeed the framework we are considering, with the additional difficulty of the absence of a mass term.
	\vskip0.2truecm
	Choquard-type equations, namely Schr\"odinger equations with a nonlocal right-hand side, appear in many physics contexts, since they originate from systems where a Schr\"odinger and a Poisson equation are coupled: those systems, indeed, model, among others, the interaction of two identically charged particles in electromagnetism, and the self-interaction of the wave function with its own gravitational field in quantum mechanics. For the physics background we refer to \cite{BF,LRZ} and to the references therein. The mathematical interest lies on the fact that the equations of the form
	\begin{equation}\label{Choq_mu_mass}\tag{$\text{Ch}$}
		-\Delta u+V(x)u=\left(\frac1{|\cdot|^\mu}\ast F(u)\right)f(u)\quad\ \mbox{in}\ \ \R^N,
	\end{equation}
	where $f$ is a subcritical or critical nonlinearity, can be treated by variational methods. Indeed, if $N\geq3$ and in the case where the potential $V>0$, one usually works in the natural Sobolev space $H^1(\R^N)$ and takes advantage of the Hardy-Littlewood-Sobolev inequality (see Lemma \ref{HLS} below) to prove that the functional associated to \eqref{Choq_mu_mass} is well-defined, see \cite{MV1,MV,CZ,CVZ}. The planar case $N=2$ is more delicate, since this setting is limiting for the Sobolev embedding, and specific techniques need to be developed, see \cite{ACTY,AFS}. Note, however, that in order to retrieve the connection with the Schr\"odinger-Poisson system, the kernel $|\cdot|^{-\mu}$ should be replaced by $-\log|\cdot|$, which is sign-changing and unbounded from above and below, and this makes the analysis even harder: we refer to \cite{CW,CT,LRTZ,CDL} and to the recent developments in \cite{BCT,CLR,CLR2}.
	\vskip0.2truecm
	Some physics model prescribe however that the potential $V$ appearing in the Schr\"odinger equation is identically zero, e.g. in the study of the Yang-Mills equation in the nonabelian gauge theory of particle physics, see \cite{Gi}. Such "zero mass case" is mathematically intriguing, since the absence of the mass implies a lack of control of the $L^2$-part of the norm in $H^1(\R^N)$. Therefore -- even if the right-hand side is just local -- one is lead to study the equation in the \textit{homogeneous Sobolev space} $D^{1,2}_0(\R^N)$, defined as the completion of $C^\infty_0(\R^N)$ with respect to the norm $\|\nabla\cdot\|_2$. In the higher dimensional case $N\geq3$, one can still work in this homogeneous space, thanks to the critical Sobolev embedding $D^{1,2}_0(\R^N)\hookrightarrow L^{2^*}(\R^N)$, see 
	e.g. \cite{BL,AP,ASM} for Schr\"odinger equations and \cite{AY} for Choquard equations. However, in the Sobolev limiting case $N=2$, where already the additional difficulty of being able to deal with exponential nonlinearities appears, the space $D^{1,2}_0(\R^2)$ is not a space of functions anymore: indeed, one cannot distinguish between $u$ and $u+c$ for all $c\in\R$, and no Sobolev embeddings can be proved in this setting. The same problem of course occurs for $D^{1,N}_0(\R^N)$. We point out nevertheless that, when dealing with Choquard equations with zero mass and logarithmic kernel, that is originating from Schr\"odinger-Poisson systems, a sort of mass term may be anyway retrieved from the nonlocal term by a careful splitting of the logarithm, and this enables one to work again in a (possibly inhomogeneous) Sobolev space, see \cite{WCR,CSTW} for the linear case $f(u)=u$ and \cite{BRT} for the delicate extension for a general class of subcritical or critical nonlinearities. This trick however does not work in the case of a polynomial kernel.
	
	In the context of Schr\"odinger equations with zero mass in $\R^N$, in the recent paper \cite{dAC} the authors managed to retrieve a good functional framework by modifying the operator, namely introducing in the standard $N$-Laplacian $\divv\left(|\nabla u|^{N-2}\nabla u\right)$ a positive radial \textit{coercive} weight function $A$, i.e. which satisfies
	\begin{enumerate}
		\item[(A)] $A:\R^+\to\R$ is continuous, $\liminf_{r\to0^+}A(r)>0$ and there exist $A_0,\ell>0$ such that $A(r)\geq A_0 r^\ell$, for all $r>0\,$,
	\end{enumerate}
	and considering the weighted operator $\divv\left(A(|x|)|\nabla u|^{N-2}\nabla u\right)$. In this case, the functional space which naturally arises is
	\begin{equation}\label{space}
		E:=\Big\{u\in L^N_{loc}(\R^N)\,\Big|\,\int_{\R^N}A(|x|)|\nabla u|^N\dd x<+\infty\Big\},
	\end{equation}
	which is a reflexive\footnote{The reflexivity of $E$ can be shown in the usual way thanks to the reflexivity of the weighted Lebesgue spaces $L^N(\R^N\!,A(\cdot)\dd x)$ for $N\geq2$ see e.g. \cite{DS}.} Banach space when endowed with the norm
	\begin{equation}\label{normE}
		\|u\|:=\left(\intN A(|x|)|\nabla u|^N\dd x\right)^\frac1N,
	\end{equation}
	see \cite[Lemma 2.1 and Corollary 1.5]{dAC}. In particular, in its radial subspace, denoted by $\Erad$, one may recover the Sobolev embeddings, which are necessary not only to enable us to accomplish our estimates, but also to show that $\Erad$ is in fact a space of functions. For $p\geq 1$, let us define first the $Q$-weighted Lebesgue space
	$$L^p_Q(\R^N):=\Big\{u\in\cM(\R^N)\,\Big|\,\intN Q(|x|)|u|^p\dd x<+\infty\Big\}\,,$$
	where $\cM(\R^N)$ stands for the set of all measurable functions on $\R^N$.
	\begin{mainth}(\cite{dAC}, Theorem 1.2)\label{Thm_cpt_emb}
		Assume (A) and 
		\begin{enumerate}
			\item[(Q)] $Q:\R^+\to\R^+$ is continuous and there exist $b_0,b>-N$ such that
			$$\limsup_{r\to0^+}\frac{Q(r)}{r^{b_0}}<+\infty\quad\mbox{and}\quad\limsup_{r\to+\infty}\frac{Q(r)}{r^b}<+\infty\,.$$
		\end{enumerate}
		Then the embedding $\Erad\hookrightarrow L^p_Q(\R^N)$ is continuous for $\gamma\leq p<+\infty$, where
		\begin{equation}\label{gamma}
			\gamma:=\max\left\{N,\frac{(b-\ell+N)(N+1)}\ell+N\right\}=\begin{cases}
				N\ \  &\mbox{if}\ \,b<\ell-N,\\
				\frac{(b-\ell+N)(N+1)}\ell+N & \mbox{if}\ \,b\geq\ell-N.\\
			\end{cases}
		\end{equation}
		Furthermore, the embedding is compact for $\gamma\leq p<+\infty$ when $b<\ell-N$, and for $\gamma<p<+\infty$ when $b\geq\ell-N$.
	\end{mainth}
	Note that assumption (Q) allows for weight functions which can be singular at the origin and vanishing at infinity, and has been used also in the study of Choquard equation with vanishing potential, see e.g. \cite{AFS}.
	
	In \cite{dAC} the authors were also able to prove that in this limiting setting a sort of Poho\v zaev-Trudinger-Moser inequality holds. The critical exponential growth is the same as in the unweighted case, namely $t\mapsto\e^{\alpha|t|^{\frac N{N-1}}}$, while the influence of the weight functions lies in the Moser exponent. Since we are considering the whole space, one need to subtract the first terms of the Taylor expansion from the exponential, by introducing the functions
	\begin{equation}\label{Phi}
		\Phi_{\alpha,j_0}(t):=\e^{\alpha|t|^{\frac N{N-1}}}-\sum_{j=0}^{j_0-1}\frac{\alpha^j}{j!}|t|^{j\frac N{N-1}}\,,
	\end{equation}
	for $\alpha>0$ and $j_0\in\N$.
	\begin{mainth}(\cite{dAC}, Theorem 1.6)\label{ThmAC_TM}
		Assume (A) and (Q) hold, and let $j_0=\inf\big\{j\in\N\,|\,j\geq\frac{\gamma(N-1)}N\big\}$. Then, for each $u\in\Erad$ and $\alpha>0$, the function $\Phi_{\alpha,j_0}(u)$ belongs to $L^1_Q(\R^N)$. Moreover, if
		\begin{equation}\label{alphaNtilde}
			0<\alpha<\widetilde\alpha_N(Q):=\alpha_N\left(1+\frac{b_0}N\right)\left(\inf_{x\in B_1(0)}A(|x|)\right)^{\frac1{N-1}},
		\end{equation}
		where $\alpha_N:=N\omega_{N-1}^{1/(N-1)}$, with $\omega_{N-1}$ denoting the measure of the unit sphere in $\R^N$, then
		$$\sup_{u\in\Erad,\,\|u\|\leq1}\intN Q(|x|)\Phi_{\alpha,j_0}(u)\dd x<+\infty\,.$$
	\end{mainth}
	With these tools available in $\Erad$, the existence for the Schr\"odinger equation with zero mass
	\begin{equation*}
		-\divv\left(A(|x|)|\nabla u|^{N-2}\nabla u\right)=Q(|x|)f(u)\quad\mbox{in}\ \,\R^N
	\end{equation*}
	was proved in \cite{dAC}, in case $f$ is a positive critical exponential nonlinearity in the sense of Theorem \ref{ThmAC_TM}, which undergoes a strong growth condition, which is effective in a neighbourhood of zero, namely 
	\begin{equation}\label{global_ass}
		F(s)\geq\lambda s^\nu\qquad\mbox{with}\ \ \nu>\gamma\ \ \mbox{and}\ \  \lambda\ \ \mbox{large enough,}
	\end{equation}
	for all $s\in\R$ and $\gamma$ defined in \eqref{gamma}. In this functional framework, a Schr\"odinger-Poisson system with zero mass, in gradient form and with critical exponential nonlinearities, was recently considered in \cite{R}. After reducing the system to the Choquard equation with logarithmic kernel
	\begin{equation}\label{Choq_log}
		-\divv\left(A(|x|)|\nabla u|^{N-2}\nabla u\right)=C_N\left(\log\frac1{|\cdot|}\ast Q(|\cdot|)F(u)\right)Q(|x|)\,f(u)\quad\ \mbox{in}\ \ \R^N,
	\end{equation}
	existence is proved using a variational approximating procedure in the spirit of \cite{LRTZ,CLR,CLR2}: in fact, the difficulties due to a sign-changing kernel which is unbounded from below and above, are overcome by means of a uniform approximation which exploits suitable kernels having a polynomial behaviour. The global condition \eqref{global_ass} was also avoided by obtaining a fine upperbound on the mountain-pass level by means of a careful analysis on Moser sequences.
	\vskip0.2truecm
	In this paper, we study \eqref{Choq_mu}, which is the counterpart of \eqref{Choq_log} where the logarithm is substituted by the polynomial kernel $|\cdot|^{-\mu}$, $\mu\in(0,N)$. On the one hand the analysis will be less involved than the one in \cite{R}, since we do not have to face the problem of a sign-changing kernel, and thus we can work directly with the equation; on the other hand, we would like not to rely on the global growth condition \eqref{global_ass}, so a fine analysis on the mountain-pass level should still be performed.
	\vskip0.2truecm
	Before stating our results, let us introduce some additional conditions on $A$ and $Q$:
	\begin{enumerate}
		\item[(Q$_\mu$)] $Q:\R^+\to\R^+$ is continuous and there exist $b_0,b>\frac\mu2-N$ such that
		$$\limsup_{r\to0^+}\frac{Q(r)}{r^{b_0}}<+\infty\quad\mbox{and}\quad\limsup_{r\to+\infty}\frac{Q(r)}{r^b}<+\infty\,;$$
		\item[(A')] there exist $r_0>0$ and $L>0$ such that $A_0(1+|x|^\ell)\leq A(|x|)\leq A_0(1+|x|^L)$ for all $x\in B_{r_0}(0)$, with $A_0,\ell$ given by (A);
		\item[(Q')] $\liminf_{r\to0^+}Q(r)/r^{b_0}=C_Q>0$.
	\end{enumerate}
	The last two conditions will be needed in estimating the mountain pass level, and can also be found in \cite{dAC,R}, while ($Q_\mu$) is the adaptation of assumption ($Q$) to the Choquard case, and is used to prove that the functional assocated to \eqref{Choq_mu} is well-defined, see Lemma \ref{J_welldefined_MP} below.
	\vskip0.2truecm
	\noindent\textbf{Notation:} With a little abuse, from now on $A(x):=A(|x|)$ and similarly $Q(x):=Q(|x|)$.
	\vskip0.2truecm
	Concerning the nonlinearity $f$, aiming at modeling both the subcritical and the critical case, we assume the following conditions:
	\begin{enumerate}
		\item[(\textit{f}$_0$)] $f\in C^1(\R)$, $f(t)>0$ for $t>0$, and $f(t)=0$ for $t\leq0\,$;
		\item[(\textit{f}$_1^{\,s}$)] $f$ is \textit{subcritical} in the sense of Trudinger-Moser, namely
		\begin{equation*}
			\lim_{t\to+\infty}\frac{f(t)}{\e^{\alpha t^{\frac N{N-1}}}}=0\quad\mbox{for all}\ \ \alpha>0\,;
		\end{equation*}
		\item[(\textit{f}$_1^{\,c}$)] $f$ is \textit{critical} in the sense of Trudinger-Moser, namely there exists $\alpha_0>0$ such that
		\begin{equation*}
			\lim_{t\to+\infty}\frac{f(t)}{\e^{\alpha t^{\frac N{N-1}}}}=\begin{cases}
				0&\quad\mbox{for}\ \ \alpha>\alpha_0\,,\\
				+\infty&\quad\mbox{for}\ \ \alpha<\alpha_0\,;\end{cases}
		\end{equation*}
		\item[(\textit{f}$_2$)] there exists $\tp>\left(1-\tfrac\mu{2N}\right)\!\gamma$ such that $f(t)=o(t^{\tp-1})$ as $t\to0^+$;
		\item[(\textit{f}$_3$)] there exist $\tau\in\left(1-\tfrac2N,1\right)$ and $C>0$ such that
		$$\tau\leq\frac{F(t)f'(t)}{(f(t))^2}\leq C\quad\mbox{for any}\ \,t>0\,;$$
		\item[(\textit{f}$_\xi$)] there exists $\xi>0$ and $\nu>\gamma$ such that 
		\begin{equation*}
			F(t)\geq\xi t^\nu\quad\mbox{for}\ \,t\in(0,1]\,;
		\end{equation*}
		\item[(\textit{f}$_4$)]there exist $t_0,M_0>0$ and $\theta\in(0,N-1]$ such that
		$$0<t^\theta F(t)\leq M_0f(t)\quad\mbox{for}\ \,t\geq t_0\,;$$
		\item[(\textit{f}$_5$)] there exists $\beta_0>0$ such that
		$$\liminf_{t\to+\infty}\,\frac{F(t)}{\e^{\alpha_0t^{\frac N{N-1}}}}\geq\beta_0>0\,.$$
	\end{enumerate}
	
	\begin{defn}[Solution of \eqref{Choq_mu}]
		We say that $u\in E$ is a \textit{weak solution of} \eqref{Choq_log} if
		\begin{equation*}
			\intN A(x)|\nabla u|^{N-2}\nabla u\nabla\varphi\dd x=\intN\!\!\left(\intN\frac{Q(y)F(u(y))}{|x-y|^\mu}\dd y\!\right)\!Q(x)f(u(x))\varphi(x)\dd x
		\end{equation*}
		for all $\varphi\in E$.
	\end{defn}
	
	\begin{thm}\label{Thm_Choq_mu}
		Let $\mu\in(0,N)$, under conditions (A), (Q$_\mu$), ($f_0$), ($f_2$), ($f_3$), assume either that
		\begin{itemize}
			\item[S)] the problem is subcritical, namely ($f_1^{\,s}$) holds, 
		\end{itemize} or that
		\begin{itemize}
			\item[C)] the problem is critical, namely ($f_1^{\,c}$) holds, and
		\begin{itemize}
			\item[$i)$] ($f_\xi$) holds with $\xi>\xi_0$ (depending on $\nu$) given in \eqref{xi}
		\end{itemize}
		or, alternatively,
		\begin{itemize}
			\item[$ii)$] (A'), (Q'), ($f_4$)-($f_5$) are fulfilled.
		\end{itemize}
		\end{itemize}
		Then \eqref{Choq_log} has a positive radially symmetric weak solution in $\Erad$.
	\end{thm}
	
	\begin{remark}
		We stress the fact that our results are new even in the planar case $N=2$. Moreover, they can be seen as an extension of the corresponding results in \cite{ACTY,AFS} to the zero mass case, of those in \cite{R} to the case of polynomial kernels, and of those in \cite{dAC} to the Choquard framework.
	\end{remark}
	
	\begin{remark}
		Since the weight $A$ is continuous and bounded below by (A), it is clear that for all $\Omega\subset\subset\R^N$ there exist constants $\underline a_0,\overline a_0>0$ such that $\underline a_0<A(x)<\overline a_0$ for all $x\in\Omega$. This implies that $E\subset D^{1,N}(\R^N)\subset W^{1,N}_{loc}(\R^N)$, where $D^{1,N}(\R^N)$ is the homogeneous Sobolev space given by \eqref{space} with $A\equiv1$, see \cite[Lemma II.6.1]{G}. Therefore, it is sufficient to prove the existence of a nonnegative solution of \eqref{Choq_log} in order to retrieve its positivity by the strong maximum principle for quasilinear equations, see \cite[Theorem 11.1]{PS}.
	\end{remark}

	\paragraph{\textbf{Notation.}} For $R>0$ and $x_0\in\R^N$ we denote by $B_R(x_0)$ the ball of radius $R$ and center $x_0$. Given a set $\Omega\subset\R^N$, its characteristic function is denoted by $\chi_\Omega$ and $\Omega^c:=\R^N\setminus\Omega$. The space of the infinitely differentiable functions which are compactly supported is $C^\infty_0(\R^N)$, while $L^p(\R^N)$ with $p\in[1,+\infty]$ is the Lebesgue space of $p$-integrable functions. The norm of $L^p(\R^N)$ is denoted by $\|\cdot\|_p$. For $q>0$ we define $\lfloor q\rfloor$ as the largest integer strictly less than $q$; if $q>1$ its conjugate H\"older exponent is $q':=\frac q{q-1}$. The symbol $\lesssim$ indicates that an inequality holds up to a multiplicative constant depending only on structural constants. Finally, $o_n(1)$ denotes a vanishing real sequence as $n\to+\infty$. Hereafter, the letter $C$ will be used to denote positive constants which are independent of relevant quantities and whose value may change from line to line.
	
	\paragraph{Overview} After the short Section \ref{Sec_2}, in which we discuss some consequences of our assumptions and state some useful results, we prove existence for the Choquard equation \eqref{Choq_log}, splitting the proof in Sections \ref{Sec_3} and \ref{Sec_4}, according to the set of assumptions considered in the Theorem \ref{Thm_Choq_mu}.
	
	\section{Preliminaries}\label{Sec_2}
	
	From now on, we set $\Phi_\alpha:=\Phi_{j_0,\alpha}$ with $j_0$ defined in Theorem \ref{ThmAC_TM}. We start by collecting some comments on our assumptions:
	\begin{remark}
		\begin{enumerate}
			\item[(i)] From ($f_0$)-($f_1^{\,c}$)-($f_2$), and \eqref{Phi}, it is easy to infer that for fixed $\alpha>\alpha_0$, $p>1$ and for any $\varepsilon>0$ one has
			\begin{equation}\label{f-C-above}
				|f(t)|\leq\varepsilon|t|^{\tp-1}+C_1(\alpha,p,\varepsilon)|t|^{p-1}\Phi_\alpha(t),\qquad t\in\R,
			\end{equation}
			for some $C_1(\alpha,p,\varepsilon)>0$, and consequently,
			\begin{equation}\label{F-C-above}
				|F(t)|\leq\varepsilon|t|^\tp+C_2(\alpha,p,\varepsilon)|t|^p\Phi_\alpha(t),\qquad t\in\R,
			\end{equation}
			for some $C_2(\alpha,p,\varepsilon)>0$. In the case ($f_1^{\,s}$) holds in place of ($f_1^{\,c}$), the inequalities \eqref{f-C-above}-\eqref{F-C-above} are valid with $\alpha>0$ arbitrary.
			\item[(ii)] Assumption ($f_3$) implies that $f$ is monotone increasing and
			\begin{equation}\label{f3_consequence}
				F(t)\leq(1-\tau)tf(t)\quad\ \mbox{for any}\ \ t\geq0\,.
			\end{equation}
			\item[(iii)] Although frequent in the literature, see e.g. \cite{AF,AFS,dAC}, assumption ($f_\xi$) is very strong, not because of the polynomial growth $t\mapsto t^\nu$ with $\nu>\gamma$, which is reasonable since it excludes just exponential decays at $0$, but mainly because of the fact that one should prescribe this behaviour in the whole range $[0,1]$ and not just asymptotically. In fact, it is not of easy verification. For instance, the easiest example $F(t)=t^\kappa$ with $\gamma<\kappa<\nu$ verifies this growth condition just in a small right neighbourhood of $0$ and not in the whole $[0,1]$. This is the reason why we are considering also an alternative proof of our main result which uses assumptions ($f_4$)-($f_5$), although the argument which exploits ($f_\xi$) is way easier.
			\item[(iv)] ($f_5$) is a condition at infinity, compatible with the critical growth ($f_1^{\,c}$) and related to the well-known de Figueiredo-Miyagaki-Ruf condition \cite{dFMR}. It is crucial in order to estimate the mountain pass level and gain compactness, see Lemma \ref{MP_level}.
			A similar condition appears also in \cite{ACTY,CT,BCT,R,BRT}, however, as in \cite{AFS,CSTW}, we do not prescribe $\beta_0$ large.
			\item[(v)] Examples of admissible nonlinearities are $F(t)=t^q\e^{t^\alpha}$ with $q>\left(1-\tfrac\mu{2N}\right)\!\gamma$ and $\alpha\in\left[0,\frac N{N-1}\right]$. The critical case corresponds to the choice of $\alpha=\frac N{N-1}\,$.
		\end{enumerate}
	\end{remark}

	\noindent The next lemma assures that the function $\Phi_\alpha$ introduced in \eqref{Phi} has the same properties of the exponential. 
	
	\begin{lem}
		For $\alpha>0$, $r>1$, and $\nu\in\R$ there holds
		\begin{equation}\label{estimate_Sani}
			(\Phi_\alpha(t))^r\leq\Phi_{\alpha r}(t)\qquad\mbox{for all}\ \,t>0
		\end{equation}
		and
		\begin{equation}\label{estimate_CFFM}
			\Phi_\alpha(\nu t)=\Phi_{\alpha\nu^{\frac N{N-1}}}(t)\qquad\mbox{for all}\ \,t>0\,.
		\end{equation}
	\end{lem}
	\begin{proof}
		For the first inequality see \cite[Lemma 2.1]{Y}; the second is just an easy calculation.	
	\end{proof}
	
	\vskip0.2truecm
	We end this section by recalling the well-known Hardy-Littlewood-Sobolev inequality, see \cite[Theorem 4.3]{LL}, which will be frequently used throughout the paper.
	\begin{lemma}(Hardy-Littlewood-Sobolev inequality)\label{HLS}
		Let $N\geq1$, $s,r>1$, and $\mu\in(0,N)$ with $\tfrac1s+\tfrac\mu N+\tfrac1r=2$. There exists a constant $C=C(N,\mu,s,r)$ such that for all $f\in L^s(\R^N)$ and $h\in L^r(\R^N)$ one has
		$$
		\int_{\R^N}\left(\frac1{|\cdot|^\mu}\ast f\right)\!h\dd x\leq C\|f\|_s\|h\|_r\,.
		$$
	\end{lemma}
	
	\section{Proof of Theorem \ref{Thm_Choq_mu}: the subcritical case and the critical case (\textit{i})}\label{Sec_3}
	We start by proving that the functional $J$, formally associated to \eqref{Choq_mu},
	\begin{equation*}
		J(u):=\frac1N\intN\!A(x)|\nabla u|^N\!\dd x-\frac12\intN\!\left(\intN\frac{Q(y)F(u(y))}{|x-y|^\mu}\dd y\right)Q(x)F(u(x))\dd x
	\end{equation*}
	is well-defined in the space $\Erad$, is $C^1$ with derivative
	$$J'(u)[\varphi]=\intN A(x)|\nabla u|^{N-2}\nabla u\nabla\varphi\dd x-\intN\!\!\left(\intN\frac{Q(y)F(u(y))}{|x-y|^\mu}\dd y\!\right)\!Q(x)f(u(x))\varphi(x)\dd x\,,$$
	and possesses a mountain-pass geometry.
	\begin{lem}\label{J_welldefined_MP}
		Under assumptions ($f_0$), ($f_2$) and either ($f_1^{\,c}$) or ($f_1^{\,s}$), the functional $J:\Erad\to\R$ is well-defined and $C^1$. If $f$ satisfies also ($f_3$), there exist constants $\rho,\eta>0$ and $e\in\Erad$ such that:
		\begin{enumerate}
			\item[{\rm(i)}] $J|_{S_\rho}\geq\eta>0$, where $S_\rho=\big\{u\in\Erad\,|\,\|u\|=\rho\big\}$;
			\item[{\rm(ii)}] $\|e\|>\rho$ and $J(e)<0\,$.
		\end{enumerate}
	\end{lem}
	\begin{proof}
		Although the proof is standard, the main tool being the Hardy-Littlewood-Sobolev inequality (Lemma \ref{HLS}), we retrace it here, in particular to show the r\^ole of assumption ($Q_\mu$).
		
		We focus on the second term of $J$, the first one being already $\|u\|^N$, see \eqref{normE}. By Lemma \ref{HLS} with $r=t=\tfrac{2N}{2N-\mu}$, the upperbound \eqref{F-C-above} and H\"older's inequality, one infers
		\begin{equation*}
			\begin{split}
				\intN\!&\left(\intN\frac{Q(y)F(u(y))}{|x-y|^\mu}\dd y\right)Q(x)F(u(x))\dd x\les\left(\intN|QF(u)|^{\frac{2N}{2N-\mu}}\right)^{\frac{2N-\mu}N}\\
				&\les\left(\intN Q^{\frac{2N}{2N-\mu}}|u|^{\frac{2N\tp}{2N-\mu}}\right)^{\frac{2N-\mu}N}\!\!+\left(\intN Q^{\frac{2N}{2N-\mu}}|u|^{\frac{2N\tp q'}{2N-\mu}}\right)^{\frac{2N-\mu}{Nq'}}\!\left(\intN Q^{\frac{2N}{2N-\mu}}|\Phi_\alpha(u)|^{\frac{2Nq}{2N-\mu}}\right)^{\frac{2N-\mu}{Nq}}\!\!,
			\end{split}
		\end{equation*}
		for $\alpha>\alpha_0$ in case ($f_1^{\,c}$) holds (resp. $\alpha>0$ if ($f_1^{\,s}$) holds). In order to use now the Sobolev embedding given by Theorem \ref{Thm_cpt_emb}, as well as to bound the exponential term by Theorem \ref{ThmAC_TM}, in both cases the weight function $\tQ:=Q^{\frac{2N}{2N-\mu}}$ must verify assumption ($Q$), and the exponent of $u$, namely $\tfrac{2N\tp}{2N-\mu}$ should be greater than $\gamma$. However, it is not difficult to show that this is the case under our assumptions ($Q_\mu$) and ($f_2$). As a result, using also \eqref{estimate_Sani}, one infers
		\begin{equation}\label{stima_convoluzione}
			\begin{split}
				\intN\!&\left(\intN\frac{Q(y)F(u(y))}{|x-y|^\mu}\dd y\right)Q(x)F(u(x))\dd x\\
				&\quad\les\|u\|^{2\tp}+\|u\|^{2p}\!\left(\intN\tQ\,\Phi_{\frac{2Nq\alpha}{2N-\mu}}(u)\right)^{\frac{2N-\mu}{Nq}}\!<+\infty\,.
			\end{split}
		\end{equation}
		This shows the well-posedness of $J$ in $\Erad$, while the regularity of $J$ follows by standard arguments. In order to show (\textit{i}), from \eqref{stima_convoluzione} and \eqref{estimate_CFFM} we deduce 
		\begin{equation*}
			\begin{split}
				J(u)\ges\|u\|^N-\|u\|^{2\tp}-\|u\|^{2p}\left(\intN\tQ\,\Phi_{\frac{2Nq\alpha}{2N-\mu}\|u\|^{\frac N{N-1}}}\left(\tfrac u{\|u\|}\right)\right)^{\frac{2N-\mu}{Nq}}.
			\end{split}
		\end{equation*}
		Therefore, in order to apply the uniform estimate of Theorem \eqref{ThmAC_TM}, one needs $\frac{2Nq\alpha}{2N-\mu}\|u\|^{\frac N{N-1}}<\widetilde\alpha_N(\tQ)$ defined in \eqref{alphaNtilde}, namely to require that $\rho<\left(\frac{2N-\mu}{2Nq\alpha}\,\widetilde\alpha_N(\tQ)\right)^{\frac{N-1}N}$. If so,
		\begin{equation*}
			\begin{split}
				J(u)\ges\|u\|^N-\|u\|^{2\tp}-\|u\|^{2p},
			\end{split}
		\end{equation*}
		which implies that $0$ is a local minimum by choosing $p$ large enough, since $2\tp>(2N-\mu)\tfrac\gamma N>N$. Let us now take $0\leq\varphi\in\Erad(\R^N)$ and define
		\begin{equation*}
			\psi(t):=\frac12\intN\left(\frac1{|\cdot|^\mu}\ast QF(t\varphi)\right)Qf(t\varphi)\dd x\,.
		\end{equation*}
		Using \eqref{f3_consequence}, it is then standard to show that $\tfrac{\psi'(t)}{\psi(t)}\geq\tfrac2{(1-\tau)t}$, which in turn implies $\psi(t)\geq\psi(1)t^{\frac2{1-\tau}}$. Hence,
		\begin{equation*}
			J(t\varphi)=\frac{t^N}N\|\varphi\|^N-\psi(t)\leq\frac{t^N}N\|\varphi\|^N-Ct^{\frac2{1-\tau}}\to-\infty,
		\end{equation*}
		since $\tau\in\left(1-\frac2N,1\right)$ by ($f_3$). It is then sufficient to take $e:=t_0\varphi$ with $t_0$ large enough, to conclude that (\textit{ii}) holds.
	\end{proof}
	
	As a consequence of this mountain-pass geometry, one infers the existence of a Cerami sequence in $\Erad$ at level 
	$$c_{mp}:=\inf_{\gamma\in\Gamma}\max_{t\in[0,1]}J(\gamma(t))\,,$$
	where
	$$\Gamma:=\left\{\gamma\in C([0,1],\Erad)\,|\,\gamma(0)=0,\gamma(1)=e\right\}\,,$$
	namely, a sequence $(u_k)_k\subset\Erad$ such that
	\begin{equation}\label{PS_sequence}
		J(u_k)\to c_{mp}\quad\text{and}\quad(1+\|u_k\|)J'(u_k)\to0\quad\text{in}\,\,\big(\Erad\big)'
	\end{equation}
	as $k\to+\infty$. In details,
	\begin{equation}\label{Cerami_J}
		J(u_k)=\frac1N\intN A(x)|\nabla u_k|^N\dd x-\frac12\intN\left(\frac1{|\cdot|^\mu}\ast QF(u_k)\right)QF(u_k)=c_{mp}+o_k(1)\,,
	\end{equation}
	and for all $\varphi\in\Erad$ one has
	\begin{equation}\label{Cerami_Jder}
		J'(u_k)[\varphi]=\intN A(x)|\nabla u_k|^{N-2}\nabla u_k\nabla\varphi-\intN\left(\frac1{|\cdot|^\mu}\ast QF(u_k)\right)Qf(u_k)\varphi=o_k(1)\|\varphi\|\,,
	\end{equation}
	from which
	\begin{equation}\label{Cerami_Jder_un}
		J'(u_k)[u_k]=\intN A(x)|\nabla u_k|^N\dd x-\intN\left(\frac1{|\cdot|^\mu}\ast QF(u_k)\right)Qf(u_k)u_k=o_k(1)\|u_k\|\,.
	\end{equation}

	\begin{lemma}\label{Lem:c-bounded}
		Assume that ($f_0$)-($f_3$) hold. Let $(u_k)_k\subset\Erad$ be a Cerami sequence of $J$ at level $c_{mp}$. Then $(u_k)_k$ is bounded in $E$ with
		\begin{equation}\label{bounds_Cerami}
			\|u_k\|^N\leq c_{mp}\left(\frac1N-\frac{1-\tau}2\right)^{-1}+o_k(1)\,,
		\end{equation}
		and there exists $u\in\Erad$ such that $u_k\rightharpoonup u$ in $\Erad\,$.
	\end{lemma}
	\begin{proof}
		By \eqref{PS_sequence} and \eqref{Cerami_Jder_un} we obtain
		\begin{equation*}
			\begin{split}
				c_{mp}+o_k(1)&=J(u_k)-\frac{1-\tau}2J'(u_k)[u_k]\\
				&=\left(\frac1N-\frac{1-\tau}2\right)\|\nabla u_k\|^N-\frac12\!\intN\!\left(\frac1{|\cdot|^\mu}\ast QF(u_k)\right)Q\left(F(u_k)-(1-\tau)f(u_k)u_k\right)\\
				&\geq\left(\frac1N-\frac{1-\tau}2\right)\|\nabla u_k\|^N
			\end{split}
		\end{equation*}
		by \eqref{f3_consequence}. The weak convergence follows since $\Erad$ is a closed subspace of a reflexive Banach space.
	\end{proof}
	To show that the limit function $u$ is indeed a weak solution of \eqref{Choq_mu}, we may prove that $u_k\to u$ in $\Erad$. This is manageable in the subcritical case. On the other hand, in the critical case, we first need a suitable uniform control on the mountain-pass level so that one can use the uniform estimate given Theorem \ref{ThmAC_TM} in order to prove the convergence of the nonlocal term in the functional, see \eqref{bounds_Cerami_good} below. Under assumption ($f_\xi$), this is relatively easy, since by taking the constant $\xi$ large enough, one can decrease the value of the mountain pass level up to the desired threshold. This is the aim of the last part of this section, which therefore contains the proof of Theorem \ref{Thm_Choq_mu} under the first set of assumptions, while we defer its proof under the more verifiable assumptions ($f_4$)-($f_5$) to Section \ref{Sec_4}.
	\vskip0.2truecm
	In the spirit of \cite{AFS} we then prove:
	\begin{lem}
		Under ($f_0$)-($f_1^{\,c}$)-($f_2$), there exists $\xi_0>0$ explicit such that, if $f$ satisfies ($f_\xi$) with $\xi>\xi_0$, then
		\begin{equation}\label{ineq_c_fxi}
			c_{mp}<\left(\frac1N-\frac{1-\tau}2\right)\left(\frac{2N-\mu}{2N}\,\frac{\widetilde\alpha_N(\tQ)}{\alpha_0}\right)^{N-1}=:c_*\,,
		\end{equation}
		from which
		\begin{equation}\label{bounds_Cerami_good}
			\|u_k\|^{\frac N{N-1}}<\frac{2N-\mu}{2N\alpha_0}\,\widetilde\alpha_N(\tQ)\,.
		\end{equation}
	\end{lem}
	\begin{proof}
		Fix a nonnegative radial function $\varphi_0\in C^\infty_0(B_1(0))$ with values in $[0,1]$ such that $\varphi_0\equiv1$ in $B_\frac12(0)$ and $|\nabla\varphi_0|\leq 2$. Then
		\begin{equation*}
			\begin{split}
				J(\varphi_0)&=\frac1N\int_{B_1(0)\setminus B_\frac12(0)}A(x)|\nabla\varphi_0|^N\dd x-\frac12\int_{B_1(0)}\left(\frac1{|\cdot|^\mu}\ast QF(\varphi_0)\right)QF(\varphi_0)\\
				&\leq\frac{\omega_N}N\left(2^N-1\right)\!\sup_{B_1(0)\setminus B_\frac12(0)}\!\!A-\,\frac{\xi^2}2\int_{B_1(0)}\left(\frac1{|\cdot|^\mu}\ast Q\varphi_0^\nu\right)Q\varphi_0^\nu\,.
			\end{split}
		\end{equation*}
		Noting that the right-hand side tends to $-\infty$ as $\xi\to+\infty$, one may take $\xi>\xi_1$, where $\xi_1$ is chosen such that
		\begin{equation*}
			\frac{\xi_1^2}2\int_{B_1(0)}\left(\frac1{|\cdot|^\mu}\ast Q\varphi_0^\nu\right)Q\varphi_0^\nu=\frac{\omega_N}N\left(2^N-1\right)\!\sup_{B_1(0)\setminus B_\frac12(0)}\!\!A\,,
		\end{equation*}
		and get $J(\varphi_0)\leq 0$. As a result, by definition of $c_{mp}$ we can estimate as follows:
		\begin{equation}\label{chain_cmp}
			\begin{split}
				c_{mp}&\leq\max_{t\in[0,1]}J(t\varphi_0)\leq\max_{t\in[0,1]}\left(t^N\!\!\int_{B_1(0)}\!A(x)\frac{|\nabla\varphi_0|^N}N\dd x-\frac{\xi^2t^{2\nu}}2\int_{B_1(0)}\!\!\left(\frac1{|\cdot|^\mu}\ast Q\varphi_0^\nu\right)Q\varphi_0^\nu\right)\\
				&\leq\frac{\xi_1^2}2\int_{B_1(0)}\!\!\left(\frac1{|\cdot|^\mu}\ast Q\varphi_0^\nu\right)Q\varphi_0^\nu\,\max_{t\in[0,1]}\left(t^N-\chi t^{2\nu}\right),
			\end{split}
		\end{equation}
		where $\chi:=\left(\tfrac\xi{\xi_1}\right)^2>1$. It is standard to prove that the map $h(t):=t^N-\chi t^{2\nu}$ achieves its maximum in $t_0:=\left(\tfrac N{2\nu\chi}\right)^{\frac1{2\nu-N}}\in(0,1)$ since $\nu>\gamma>N$. Hence, inserting $h(t_0)$ in \eqref{chain_cmp}, one gets
		\begin{equation}\label{cmp_c*}
			c_{mp}\leq
			\frac{\xi_1^{\frac{4\nu}{2\nu-N}}}{\xi^{\frac{2N}{2\nu-N}}}\left(\frac N{2\nu}\right)^{\frac N{2\nu-N}}\frac{2\nu-N}{4\nu}\int_{B_1(0)}\!\!\left(\frac1{|\cdot|^\mu}\ast Q\varphi_0^\nu\right)Q\varphi_0^\nu=:c_0(\nu.N,\xi_1,Q,\varphi_0)\xi^{-\frac{2N}{2\nu-N}}.
		\end{equation}
		To show \eqref{ineq_c_fxi} we then need to choose $\xi$ so that the right-hand side is below the threshold $c_*$, namely
		\begin{equation}\label{xi}
			\xi>\xi_0:=\max\{\xi_1,\xi_*\}\,,
		\end{equation}
		where $\xi_*$ satisfies the equality in \eqref{cmp_c*}. At this point, combining the uniform bounds in \eqref{bounds_Cerami} and \eqref{ineq_c_fxi}, it is immediate to infer a nice uniform control on the norm of $(u_k)_k$ given by \eqref{bounds_Cerami_good}.
	\end{proof}
	
	We are now ready to prove Theorem \ref{Thm_Choq_mu} under assumptions ($f_0$)-($f_3$) and ($f_\xi$) with $\xi>\xi_0$ defined in \eqref{xi}.
	\begin{proof}[Proof of Theorem \ref{Thm_Choq_mu} (S)-(C-i)]
		We aim at proving that 
		\begin{equation}\label{aim_Thm1}
			T(u_k):=\intN\left(\frac1{|\cdot|^\mu}\ast QF(u_k)\right)Qf(u_k)(u_k-u)\to0
		\end{equation}
		as $n\to+\infty$. Indeed, if so, by \eqref{Cerami_Jder} with $\varphi=u$ and \eqref{Cerami_Jder_un}, one would infer
		\begin{equation*}
			\intN A(x)|\nabla u_k|^{N-2}\nabla u_k\nabla(u_k-u)\dd x\to0,
		\end{equation*}
		which, combined with
		\begin{equation*}
			\intN A(x)|\nabla u|^{N-2}\nabla u\nabla(u_k-u)\dd x\to0,
		\end{equation*}
		by weak convergence, would guarantee that $u_k\to u$ strongly in $E$ by means of the simple inequality (see \cite[inequality (2.2)]{Simon})
		$$(|y_1|^{N-2}y_1-|y_2|^{N-2}y_2)(y_1-y_2)\geq C(N)|y_1-y_2|^N\quad\mbox{for all}\ \,y_1,\, y_2\in\R^N.$$
		Since the functional is $C^1$, the fact that $u$ is a weak solution of \eqref{Choq_mu} directly follows.
		
		Hence, we are lead to show \eqref{aim_Thm1}. Estimating by the Hardy-Littlewood-Sobolev inequality, we obtain
		\begin{equation}\label{aim_Thm1_1}
			|T(u_k)|\leq\|QF(u_k)\|_{\frac{2N}{2N-\mu}}\|Qf(u_k)(u_k-u)\|_{\frac{2N}{2N-\mu}}\,,
		\end{equation}
		and we are proving that the first term is uniformly bounded, while the second converges to $0$. Indeed, similarly to \eqref{stima_convoluzione}, we have
		\begin{equation}\label{QF}
			\|QF(u_k)\|_{\frac{2N}{2N-\mu}}\les\|u_k\|^\tp+\|u_k\|^p\left(\intN Q^\frac{2N}{2N-\mu}\,\Phi_{\frac{2Nq\alpha}{2N-\mu}\|u_k\|^{\frac N{N-1}}}\Big(\tfrac{u_k}{\|u_k\|}\Big)\right)^\frac{2N-\mu}{2N}\!.
		\end{equation}
		If ($f_1^{\,s}$) holds, since $\|u_k\|$ is uniformly bounded by Lemma \ref{Lem:c-bounded}, then
		\begin{equation}\label{unif_bound_exponent}
			\frac{2Nq\alpha}{2N-\mu}\|u_k\|^{\frac N{N-1}}<\widetilde{\alpha}_N(\tQ)
		\end{equation}
		follows by taking a sufficiently small $\alpha>0$. On the other hand, in the critical case ($f_1^{\,c}$), by \eqref{bounds_Cerami_good} one may take $q>1$ close to $1$ and $\alpha>\alpha_0$ close to $\alpha_0$, so that \eqref{unif_bound_exponent} holds. In both cases the last term in \eqref{QF} is then bounded uniformly in $k$. As a result,
		\begin{equation}\label{QF_bdd}
			\|QF(u_k)\|_{\frac{2N}{2N-\mu}}\leq C
		\end{equation}
		by Lemma \ref{Lem:c-bounded}. Similarly, recalling the notation $\tQ:=Q^\frac{2N}{2N-\mu}$, by \eqref{f-C-above} and the H\"older inequality with conjugate exponents $\tp,\tp'=\tfrac\tp{\tp-1}$ for the first term, and $r,r'$ and $\nu,\nu'$ for the second, we get
		\begin{equation*}
			\begin{split}
				\|Q&f(u_k)(u_k-u)\|_{\frac{2N}{2N-\mu}}^{\frac{2N}{2N-\mu}}\les\left(\intN\tQ|u_k|^{\frac{2N\tp}{2N-\mu}}\right)^\frac{\tp-1}{\tp}\left(\intN\tQ|u_k-u|^\frac{2N\tp}{2N-\mu}\right)^\frac1\tp\\
				&+\left(\intN\tQ|u_k|^{(p-1)\frac{2Nr'}{2N-\mu}}\right)^\frac1{r'}\!\left(\intN\tQ|u_k-u|^\frac{2Nr\nu'}{2N-\mu}\right)^\frac1{r\nu'}\!\left(\intN\tQ\,\Phi_{\frac{2Nr\nu\alpha}{2N-\mu}\|u_k\|^{\frac N{N-1}}}\Big(\tfrac{u_k}{\|u_k\|}\Big)\right)^\frac1{r\nu}\!.
			\end{split}
		\end{equation*}
		As before, in the subcritical case, again a choice of $\alpha$ small enough is sufficient to control the exponential term, while in the critical case one needs to choose $r,\nu>1$ close to $1$ and $\alpha>\alpha_0$ close to $\alpha_0$, and consider the upperbound \eqref{bounds_Cerami_good}; in both cases we may show the boundedness of the exponential term; up to a smaller $\nu$ and a bigger $p$, one also has $\frac{2Nr\nu'}{2N-\mu}>\gamma$ and $(p-1)\frac{2Nr'}{2N-\mu}>\gamma$. Hence,
		\begin{equation}\label{Qf(u_k-u)to0}
			\begin{split}
				\|Q&f(u_k)(u_k-u)\|_{\frac{2N}{2N-\mu}}\les\|u_k\|^{\tp-1}\|u_k-u\|_{L^\frac{2N\tp}{2N-\mu}_\tQ}+\|u_k\|^{p-1}\|u_k-u\|_{L^\frac{2Nr\nu'}{2N-\mu}_\tQ}\to0
			\end{split}
		\end{equation}
		by Lemma \ref{Lem:c-bounded} and the compact embedding given by Theorem \ref{Thm_cpt_emb}. Combining \eqref{QF_bdd} and \eqref{Qf(u_k-u)to0} with \eqref{aim_Thm1_1}, \eqref{aim_Thm1} holds, and the strong convergence $u_k\to u$ follows, which proves that $u$ is a weak solution of \eqref{Choq_mu}.
	\end{proof}

	\section{Proof of Theorem \ref{Thm_Choq_mu}: the critical case (\textit{ii})}\label{Sec_4}
	
	As we mentioned in the introduction, the global growth assumption ($f_\xi$), introduced in the critical case, is in fact not of easy verification. In this section we prove the existence of a weak solution of \eqref{Choq_mu} in the critical case by using ($f_4$)-($f_5$) instead of ($f_\xi$); however, we will need also some control from below of the weight functions $A$ and $Q$ as in ($A'$)-($Q'$). The argument, inspired by \cite{ACTY,AFS} exploits the concentration behaviour of the Moser sequences to infer a suitable uniform bound for $\|u_k\|$, which turns out to depend on all coefficients and parameters in the equation. This will allow us to show the existence of a nontrivial solution for \eqref{Choq_mu}.
	\vskip0.2truecm
	Let us introduce the Moser sequence as
	\begin{equation*}
		\tw_n(x):=\begin{cases}
			(\log n)^{1-\frac1N}\ \  &\mbox{if}\ \,0\leq|x|\leq\tfrac\rho n\,,\\
			\tfrac{\log\tfrac\rho{|x|}}{(\log n)^\frac1N} & \mbox{if}\ \,\tfrac\rho n<|x|<\rho\,,\\
			0 & \mbox{if}\ \,|x|\geq\rho\,,\\
		\end{cases}
	\end{equation*}
	where $\rho\leq r_0$ is given by (A'). Using (A') we estimate from below its norm in $E$ as
	\begin{equation*}
		\begin{split}
			\intN\!A(x)|\nabla\tw_n|^N\dd x&=\frac{\omega_{N-1}}{\log n}\int_{\tfrac\rho n}^n\frac{A(r)}r\dd r\geq\frac{\omega_{N-1}A_0}{\log n}\int_{\tfrac\rho n}^n\frac{1+r^\ell}r\dd r\\
			&=\omega_{N-1}A_0\left(1+\frac{\rho^\ell}{\ell\log n}+o\left(\frac1{\log n}\right)\right),
		\end{split}
	\end{equation*}
	and analogously from above, hence we can state that
	\begin{equation}\label{Moser_norm}
		\|\tw_n\|^N=\omega_{N-1}A_0(1+\delta_n), \quad\ \ \mbox{with}\quad\ \frac{\rho^\ell/\ell}{\log n}+o\left(\frac1{\log n}\right)\leq\delta_n\leq\frac{\rho^L/L}{\log n}+o\left(\frac1{\log n}\right).
	\end{equation}
	Hence defining
	$$w_n:=\frac{\tw_n}{\left(\omega_{N-1}A_0(1+\delta_n)\right)^\frac1N}\,,$$
	one has $\|w_n\|=1$ for all $n\in\N$.
	
	\begin{lem}\label{MP_level}
		Under (A)-(A'), (Q$_\mu$)-(Q'), ($f_0$)-($f_3$), and ($f_4$)-($f_5$), one has
		\begin{equation}\label{MP_est}
			c_{mp}<\frac{\omega_{N-1}A_0}N\left(\frac{2\tb_0+2N-\mu}{2\alpha_0}\right)^{N-1}.
		\end{equation}
	\end{lem}
	\begin{proof}
		We aim at showing that there exist a suitable $\CyrB>0$ (to be chosen later) and $n_0\in\N$ such that
		\begin{equation}\label{Be}
			\max_{t\geq0}J(tw_{n_0})<\CyrB\,.
		\end{equation}
		Suppose by contradiction that \eqref{Be} does not hold. This means that for all $n\in\N$ there exists $t_n>0$ such that
		\begin{equation*}
			J(t_nw_n)=\max_{t\geq0}J(tw_n)\geq\CyrB\,.
		\end{equation*}
		Since the convolution term is positive and $\|w_n\|=1$ for all $n\in\N$, this implies
		\begin{equation}\label{MP_est_1a}
			t_n^N\geq N\CyrB\,.
		\end{equation}
		On the other hand, one may suppose that $t_n$ is chosen such that $J(t_nw_n)=\max\{J(tw_n)\,|\,t>0\}$ by the geometry of the functional on radial functions with compact support given by in Lemma \eqref{J_welldefined_MP}. Hence $\tfrac\dd{\dd t}\big|_{t=t_n}J(tw_n)=0$, from which
		\begin{equation}\label{MP_est_1b}
			t_n^N=\intN\left(\frac1{|\cdot|^\mu}\ast QF(t_nw_n)\right)Qf(t_nw_n)t_nw_n\,.
		\end{equation}
		Using assumptions ($f_4$)-($f_5$), for all $\varepsilon>0$ fixed there exists $t_\varepsilon>0$ such that for $t>\max\{t_0,t_\varepsilon\}$ one has
		\begin{equation}\label{f4f5}
			tf(t)F(t)\geq\frac{t^{\theta+1}}{M_0}(F(t))^2\geq\frac{\beta_0^2-\varepsilon}{M_0}t^{\theta+1}\e^{2\alpha_0t^{\frac N{N-1}}}.
		\end{equation}
		Hence, recalling that $w_n$ is constant in $B_{\frac\rho n}(0)$, we can estimate the right-hand side of \eqref{MP_est_1b} from below by \eqref{f4f5} as
		\begin{equation}\label{MP_est_1b_continuazione}
			\begin{split}
				t_n^N&\geq\int_{B_{\frac\rho n}(0)}\left(\int_{B_{\frac\rho n}(0)}\frac{Q(y)F(t_nw_n(y))}{|x-y|^\mu}\dd y\right)Q(x)f(t_nw_n(x))t_nw_n(x)\dd x\\
				&\geq\frac{(\beta_0^2-\varepsilon)t_n^{\theta+1}\left(\log n\right)^{\left(1-\frac1N\right)(\theta+1)}}{M_0\left((1+\delta_n)A_0\omega_{N-1}\right)^\frac{\theta+1}N}\,\e^{\frac{2\alpha_0t_n^{\frac N{N-1}}\log n}{\left((1+\delta_n)A_0\omega_{N-1}\right)^\frac1{N-1}}}\int_{B_{\frac\rho n}(0)}\!\int_{B_{\frac\rho n}(0)}\!\frac{Q(x)Q(y)}{|x-y|^\mu}\dd x\dd y\,.
			\end{split}
		\end{equation}
		By (Q') we can estimate from below $Q(r)>cr^{b_0}$ in $B_{\frac\rho n}(0)$ for $n$ large enough; hence, using the simple estimate $\frac1{|x-y|^\mu}\geq\Big(\tfrac n{2\rho}\Big)^\mu$ for all $x,y\in B_{\frac\rho n}(0)$, we obtain
		\begin{equation*}
			\begin{split}
				\int_{B_{\frac\rho n}(0)}\!&\int_{B_{\frac\rho n}(0)}\!\frac{Q(x)Q(y)}{|x-y|^\mu}\dd x\dd y\geq c^2\left(\frac n{2\rho}\right)^\mu\left(\int_{B_{\frac\rho n}(0)}|x|^{b_0}\dd x\right)^2\\
				&=c^2\left(\frac n{2\rho}\right)^\mu\omega_{N-1}^2\left(\int_0^{\frac\rho n}r^{b_0+N-1}\dd r\right)^2=\frac{c^2\omega_{N-1}^2}{2^\mu(b_0+N)^2}\left(\frac\rho n\right)^{2b_0+2N-\mu}.
			\end{split}
		\end{equation*}
		Hence, from \eqref{MP_est_1b_continuazione} one infers
		\begin{equation}\label{t_n_from_below}
			\begin{split}
				t_n^{N-\theta-1}\geq K\exp\!\left\{\!\left(\frac{2\alpha_0\,t_n^{\frac N{N-1}}}{\left(A_0\omega_{N-1}(1+\delta_n)\right)^{\frac1{N-1}}}-(2b_0+2N-\mu)\!\right)\!\log n+\frac{N-1}N(\theta+1)\log\log n\right\},
			\end{split}
		\end{equation}
		where the constant $K$ is defined as
		\begin{equation*}
			K:=\frac{(\beta_0^2-\varepsilon)c^2\omega_{N-1}^2\rho^{2b_0+2N-\mu}}{M_02^\mu(b_0+N)^2\left(A_0\omega_{N-1}(1+\delta_n)\right)^{\frac{\theta+1}N}}\,.
		\end{equation*}
		Applying the $\log$ on both sides of \eqref{t_n_from_below} yields
		\begin{equation}\label{t_n_from_below_2}
			\begin{split}
				(N-1-\theta)\frac{(N-1)}Nt_n^\frac N{N-1}&\geq(N-1-\theta)\log(t_n)\geq\log K+\frac{N-1}N(\theta+1)\log\log n\\
				&\quad+\left(\frac{2\alpha_0\,t_n^\frac N{N-1}}{\left(A_0\omega_{N-1}(1+\delta_n)\right)^{\frac1{N-1}}}-(2b_0+2N-\mu)\right)\log n\,.
			\end{split}
		\end{equation}
		Dividing by $t_n^{\frac N{N-1}}$, we obtain
		\begin{equation*}
			\begin{split}
				(N-1-\theta)\frac{(N-1)}N&\geq\left(\frac{2\alpha_0}{\left(A_0\omega_{N-1}(1+\delta_n)\right)^{\frac1{N-1}}}-\frac{2b_0+2N-\mu}{t_n^\frac N{N-1}}\right)\log n\,.
			\end{split}
		\end{equation*}
		If $t_n\to+\infty$, then one would get a contradiction for large $n$, since $\theta\in(0,N-1]$. Same, if the factor in front of $\log n$ is positive. Hence we infer that $(t_n)_n$ is bounded with
		\begin{equation}\label{MP_est_fact2}
			t_n^N\leq A_0\omega_{N-1}(1+\delta_n)\left(\frac{2b_0+2N-\mu}{2\alpha_0}\right)^{N-1}\!.
		\end{equation}
		Comparing \eqref{MP_est_1a} and \eqref{MP_est_fact2}, and since $\delta_n=o_n(1)$ as $n\to+\infty$, we see that choosing
		\begin{equation}\label{Be_def}
			\CyrB:=\frac{A_0\omega_{N-1}}N\left(\frac{2b_0+2N-\mu}{2\alpha_0}\right)^{N-1}\!,
		\end{equation}
		one reaches the claimed contradiction, namely one gets
		\begin{equation*}
			\exists\lim_{n\to+\infty}t_n=A_0\omega_{N-1}\left(\frac{2b_0+2N-\mu}{2\alpha_0}\right)^{N-1}\!.
		\end{equation*}
		Now, combining \eqref{MP_est_1a}, \eqref{Be_def}, and \eqref{MP_est_fact2}, from \eqref{t_n_from_below_2} we deduce
		\begin{equation*}
			\begin{split}
				C&\geq \left(\frac{2\alpha_0\,t_n^\frac N{N-1}}{\left(A_0\omega_{N-1}(1+\delta_n)\right)^{\frac1{N-1}}}-(2b_0+2N-\mu)\right)\log n+\frac{N-1}N(\theta+1)\log\log n\\
				&\geq(2b_0+2N-\mu)\left(\frac1{(1+\delta_n)^\frac1{N-1}}-1\right)\log n+\frac{N-1}N(\theta+1)\log\log n\\
				&\geq(2b_0+2N-\mu)\left(\frac{-\delta_n}{N-1}+o(\delta_n)\right)\log n+\frac{N-1}N(\theta+1)\log\log n\\
				&=o_n(1)+\frac{N-1}N(\theta+1)\log\log n\,,
			\end{split}
		\end{equation*}
		recalling \eqref{Moser_norm}, which is again a contradiction. Therefore, \eqref{Be} with \eqref{Be_def} must hold true, which readily implies \eqref{MP_est}.
	\end{proof}
	
	With the fine upperbound of the mountain-pass level given by Lemma \ref{MP_level} we are in the position to prove the existence of a nontrivial weak solution of \eqref{Choq_mu}. The argument follows the line of \cite{ACTY}, see also \cite{AFS}, and we only sketch it, but paying attention to the more delicate points.
	
	\begin{proof}[Proof of Theorem \ref{Thm_Choq_mu}(C-ii)]
		First, we prove that
		\begin{equation}\label{conv_Ff}
			\left(\frac1{|\cdot|^\mu}\ast QF(u_k)\right)Qf(u_k)\varphi\to\left(\frac1{|\cdot|^\mu}\ast QF(u)\right)Qf(u)\varphi\qquad\mbox{in}\ \,L^1(\R^N)
		\end{equation}
		for all test functions $\varphi$, where $u$ is the limit point of the Cerami sequence $(u_k)_k$. For such $\varphi$, it is easy to prove that $w_k:=\tfrac{\varphi}{1+u_k}\in\Erad$. Indeed,
		\begin{equation*}
			\begin{split}
				\|w_n\|^N&\leq\intN A(x)\left(\frac{|\nabla\varphi|^N}{(1+u_k)^N}+\frac{|\varphi|^N|\nabla u_k|^N}{(1+u_k)^{2N}}\right)\!\dd x\\
				&\leq\intN A(x)|\nabla\varphi|^N\dd x+C(\varphi)\intN A(x)|\nabla u_k|^N\dd x\les\|\varphi\|^N+\|u_k\|^N\leq C\,.
			\end{split}
		\end{equation*}
		by Lemma \ref{Lem:c-bounded}. This implies that one may test \eqref{Cerami_Jder} with $w_k$ and find
		\begin{equation}\label{primopezzo}
			\begin{split}
				\intOmega\bigg(\frac1{|\cdot|^\mu}&\ast QF(u_k)\bigg)Qf(u_k)\frac{\varphi}{1+u_k}\dd x=\intN A(x)|\nabla u_k|^{N-2}\nabla u_k\nabla w_k\dd x+o_k(1)\|w_k\|\\
				&\leq\intN A(x)|\nabla u_k|^N|\varphi|\dd x+\intN A(x)|\nabla u_k|^{N-1}\frac{|\nabla\varphi|}{1+u_k}\dd x+o_k(1)\left(\|\varphi\|+\|u_k\|\right)\\
				&\leq2\|u_k\|^N+\|\varphi\|^N+o_k(1)\leq C\,,
			\end{split}
		\end{equation}
		since $u_k\geq0$ in the second integral, and having used the H\"older inequality there. Let $\Omega\subset\subset\R^N$ and $\varphi\geq0$ be a test function such that $\varphi\equiv1$ on $\Omega$. Then
		\begin{equation*}
			\begin{split}
				&\intOmega\left(\frac1{|\cdot|^\mu}\ast QF(u_k)\right)Qf(u_k)\dd x\\
				&\quad\leq2\int_{\{u_k\leq1\}\cap\Omega}\!\left(\frac1{|\cdot|^\mu}\ast QF(u_k)\right)\frac{Qf(u_k)}{1+u_k}+\int_{\{u_k\geq1\}\cap\Omega}\!\left(\frac1{|\cdot|^\mu}\ast QF(u_k)\right)Qf(u_k)u_k\\
				&\quad\leq\intOmega\!\bigg(\frac1{|\cdot|^\mu}\ast QF(u_k)\bigg)Qf(u_k)\frac{\varphi}{1+u_k}+\intN\!\left(\frac1{|\cdot|^\mu}\ast QF(u_k)\right)Qf(u_k)u_k\leq C\,,
			\end{split}
		\end{equation*}
		thanks to \eqref{primopezzo} and \eqref{Cerami_Jder_un}-\eqref{bounds_Cerami}. As a result, the measure $\nu_n$ defined by
		$$\nu_n(\Omega):=\intOmega\left(\frac1{|\cdot|^\mu}\ast QF(u_k)\right)Qf(u_k)\dd x$$
		has uniformly bounded total variation, hence there exists a measure $\nu$ such that, up to a subsequence, $\nu_n\stackrel{*}{\rightharpoonup}\nu$, namely
		$$\intOmega\left(\frac1{|\cdot|^\mu}\ast QF(u_k)\right)Qf(u_k)\varphi\dd x\to\intOmega\varphi\dd\nu$$
		for all $\varphi\in C^\infty_0(\Omega)$. As in \cite[Lemma 2.4]{ACTY} we may then conclude that $\nu$ is absolutely continuous with respect to the Lebesgue measure and it can be identified as $\nu=\left(\frac1{|\cdot|^\mu}\ast QF(u)\right)Qf(u)\dd x$, which proves \eqref{conv_Ff}.
		
		Combining \eqref{conv_Ff} with the weak convergence $u_k\rightharpoonup u$ in $E$, we infer that $u$ is a weak solution of \eqref{Choq_mu}. We need now to prove that $u\not\equiv0$. To this aim, we first show that 
		\begin{equation}\label{conv_FF}
			\intN\!\left(\frac1{|\cdot|^\mu}\ast QF(u_k)\right)QF(u_k)\to\intN\!\left(\frac1{|\cdot|^\mu}\ast QF(u)\right)QF(u)\,.
		\end{equation}
		Reasoning as in \cite[Lemma 2.4]{ACTY}, thanks to ($f_4$) it is possible to reduce the proof of \eqref{conv_FF} to 
		\begin{multline}\label{conv_FF_bdd}
			\int_{\{u_k\leq M\}}\!\!\left(\int_{\{u_k\leq K\}}\frac{Q(y)F(u_k(y))}{|x-y|^\mu}\dd y\right)\!Q(x)F(u_k(x))\dd x\\
			\to\int_{\{u\leq M\}}\!\!\left(\int_{\{u\leq K\}}\frac{Q(y)F(u(y))}{|x-y|^\mu}\dd y\right)\!Q(x)F(u(x))\dd x\,,
		\end{multline}
		for all $M,K>0$ large enough. However, if $u_k$ is pointwisely bounded, by ($f_2$) one deduces $F(u_k)\leq C_{M,K}|u_k|^\tp$, therefore,
		\begin{equation}\label{conv_FF_bdd_proof}
			\begin{split}
				\int_{\{u_k\leq M\}}\!\!\left(\int_{\{u_k\leq K\}}\frac{Q(y)F(u_k(y))}{|x-y|^\mu}\dd y\right)\!Q(x)F(u_k(x))\dd x&\les\|Q|u_k|^\tp\|_{\frac{2N}{2N-\mu}}^2\to\|Q|u|^\tp\|_{\frac{2N}{2N-\mu}}^2
			\end{split}
		\end{equation}
		by the strong convergence given by Theorem \ref{Thm_cpt_emb}. Hence, by the inverse of the dominated convergence theorem \cite[Theorem 1.2.7]{BS}, the left-hand side of \eqref{conv_FF_bdd_proof} is uniformly bounded and we can use the dominated convergence theorem to prove \eqref{conv_FF_bdd}, and in turn \eqref{conv_FF}.
		
		Assuming by contradiction $u\equiv0$, then combining \eqref{conv_FF}, $F(0)=0$, and \eqref{Cerami_J} one infers
		\begin{equation}\label{cmp_assurdo}
			\begin{split}
				c_{mp}&=J(u_k)+o_k(1)\\
				&=\frac{\|u_k\|^N}N+\frac12\intN\!\left(\frac1{|\cdot|^\mu}\ast QF(u_k)\right)QF(u_k)+o_k(1)=\frac{\|u_k\|^N}N+o_k(1)\,,
			\end{split}
		\end{equation}
		from which, by Lemma \ref{MP_level},
		\begin{equation}\label{stima_norma_via_cmpMoser}
			\begin{split}
				\frac{2N\alpha_0}{2N-\mu}\|u_k\|^{\frac N{N-1}}&=\frac{2N\alpha_0}{2N-\mu}(Nc_{mp})^\frac1{N-1}+o_k(1)\\
				&<\frac{2N\alpha_0}{2N-\mu}(\omega_{N-1}A_0)^\frac1{N-1}\,\frac{2b_0+2N-\mu}{2\alpha_0}\\
				&=N(\omega_{N-1}A_0)^\frac1{N-1}\left(1+\frac{2b_0}{2N-\mu}\right).
			\end{split}
		\end{equation}
		By \eqref{Cerami_Jder_un} and the Hardy-Littlewood inequality we have
		\begin{equation}\label{norma_uk_assurdo}
			\|u_k\|^N+o_k(1)\leq\|QF(u_k)\|_{\frac{2N}{2N-\mu}}\|Qf(u_k)u_k\|_{\frac{2N}{2N-\mu}}
		\end{equation}
		and we estimate the two terms as in \eqref{QF} thanks to \eqref{F-C-above} and \eqref{f-C-above}, respectively. The exponential term is then uniformly bounded by Theorem \ref{ThmAC_TM} by \eqref{stima_norma_via_cmpMoser}, since
		$$\widetilde{\alpha}_N(\tQ)=N(\omega_{N-1}A_0)^\frac1{N-1}\left(1+\frac1Nb_0\frac{2N}{2N-\mu}\right).$$
		 Since $u_k\to0$ in $L_\tQ^t(\R^N)$ for $t>\gamma$, from \eqref{norma_uk_assurdo} we conclude that $\|u_k\|\to0$, which lead us to a contradiction with \eqref{cmp_assurdo}. We can thus conclude that the weak solution $u$ is nontrivial.
	\end{proof}

\vskip0.4truecm
\paragraph{Acknowledgements:} The Author is member of \textit{Gruppo Nazionale per l'Analisi Matematica, la Probabilità e le loro Applicazioni} (GNAMPA) of the \textit{Instituto Nazionale di Alta Matematica} (INdAM), and was partially supported by INdAM-GNAMPA Project 2023 titled \textit{Interplay between parabolic and elliptic PDEs} (codice CUP E53C2200l93000l).

\end{document}